\documentstyle[11pt]{article}
\newtheorem{prop}{Proposition}[section]
\newtheorem{th}[prop]{Theorem}
\newtheorem{cor}[prop]{Corollary}

\newcommand{\bsquare}{\hbox{\rule{6pt}{6pt}}}
\newcommand{\proof}[1]{\noindent{\bf Proof}\hspace{0.3cm}{#1}\hfill\bsquare 
\vspace{0.5cm}\par}

\date{}

\begin{document}

%\title{A note on the supersingular K3 surface with Artin invariant 1 in characteristic 2}
\title{A note on a supersingular K3 surface in characteristic 2}

\author{Toshiyuki Katsura\thanks{Partially supported by JSPS Grant-in-Aid (S), No 19104001} and Shigeyuki Kond\=o\thanks{Partially supported by JSPS Grant-in-Aid (S), No 22224001}}

\maketitle

\rightline{\it Dedicated to Gerard van der Geer on the occasion of  his 60th birthday}

\begin{abstract}
We construct, on a supersingular K3 surface with Artin invariant 1 in characteristic 2,
a set of 21 disjoint smooth rational curves and another
set of 21 disjoint smooth rational curves such that each curve in one set intersects
exactly 5 curves from the other set with multiplicity 1 by using the structure of 
a generalized Kummer surface. As a corollary we have a concrete construction of 
a K3 surface with 21 rational double points of type $A_{1}$ in characteristic 2.

\end{abstract}

%contents

\section{Introduction}
Let $k$ be an algebraically closed field of characteristic $p$.
Let $A$ be an abelian surface over $k$ and let $\iota$ be the inversion
of $A$. In case $p \geq 3$ or 0, the minimal resolution
of the quotient surface $A/\langle \iota \rangle$ is a K3 surface.
It is called a Kummer surface and is denoted by ${\rm Km}(A)$. 
If $p\geq 3$ and the abelian surface $A$
is superspecial, i.e., $A$ is a product of two supersingular elliptic 
curves, then ${\rm Km}(A)$ is a supersingular K3 surface 
with Artin invariant 1, which is unique up to isomorphism.  However, in case $p = 2$
the situation is different.  Namely, if $A$ is a supersingular
abelian surface in characteristic 2, then the quotient surface
$A/\langle \iota \rangle$ is birationally equivalent to a rational surface,
and so it is not birationally equivalent to a K3 surface (Shioda \cite{Sh} and Katsura \cite{Ka}).
Therefore, we need to consider a different way to construct
a supersingular K3 surface with Artin invariant 1 in characteristic 2,
and Schr\"oer proved a certain generalized Kummer surface is isomorphic to
such a surface (cf. Schr\"oer \cite{Sch}).

On the other hand, 
Dolgachev and the second author gave several ways of constructions of the supersingular 
K3 surface $S$ with Artin invariant 1 in characteristic 2 (Dolgachev and Kondo \cite{DK}).
They showed that $S$ contains a set $\cal{A}$ of 21 disjoint smooth rational curves and another
set $\cal{B}$ of 21 disjoint smooth rational curves such that each curve in one set intersects
exactly 5 curves from the other set with multiplicity 1.  We call the configuration of these 42 curves a $(21)_5$-configuration. 
These 42 curves can be obtained by the following two ways.  
Firstly one can contract 21 curves in one family, for example in $\cal A$, and obtains a surface 
$\bar{S}$ with 21 nodes.  Then $\bar{S}$ is 
the inseparable double cover of ${\bf P}^2$ with branch divisor 
$x_0x_1x_2(x_0^3+x_1^3+x_2^3)= 0$.  The set $\cal{A}$ consists of exceptional curves of 
the minimal resolution of singularities of $\bar{S}$ and
$\cal{B}$ is the set of curves which are preimages of the lines on ${\bf P}^2$ defined over
${\bf F}_{4}$.  Secondly $S$ is isomorphic to the surface in ${\bf P}^2\times {\bf P}^2$ defined by
the equations $x_0y_0^2+x_1y_1^2+x_2y_2^2 = 0, \  x_0^2y_0+x_1^2y_1+x_2^2y_2 = 0.$
The 42 lines 
$$(x_0,x_1,x_2) = (a_0,a_1,a_2)\in {\bf P}^2({\bf F}_{4}), \ a_0^2y_0+a_1^2y_1+a_2^2y_2 = 0,$$
$$(y_0,y_1,y_2) = (a_0,a_1,a_2)\in {\bf P}^2({\bf F}_{4}), \ x_0a_0^2+x_1a_1^2+x_2a_2^2 = 0$$
give the two sets $\cal{A}$ and $\cal{B}$.  

The main purpose in this paper is to give another construction
of 42 curves on $S$ as above by using
a generalized Kummer surface.
We shall show that $S$ is the minimal resolution of the quotient of $E\times E$ 
by an automorphism of order 3 where $E$ is the supersingular elliptic curve.
The quotient surface has 9 rational double points of type $A_2$.  Thus $S$ contains
18 smooth rational curves obtained by resolution of singularities.  On the other hand 
we shall show that there exist 24 elliptic curves on $E\times E$ whose images, together
with 18 exceptional curves, form a $(21)_5$-configuration (Theorem 4.1).

This and the above two constructions of $S$ give an analogue
of $(16)_6$-configuration on Kummer surfaces associated with a curve of genus
two.  Let $C$ be a smooth curve of genus two defined over $k$ with $p\geq 3$ or $0$.
Let $J(C)$ be the Jacobian of $C$.  Then
the quotient surface $J(C)/\langle \iota \rangle$ can be embedded into ${\bf P}^{3}$
as a quartic surface with $16$ nodes.  There are sixteen planes in ${\bf P}^{3}$ each of which
touches $J(C)/\langle \iota \rangle$ along a conic.  These conics are the 
images of the theta divisor on $A$ and its translations
by 2-torsion points. Thus ${\rm Km}(J(C))$ contains two sets of 16 disjoint smooth rational curves.
Each member in one set meets exactly six members in another family. 
The configuration of these thirty-two curves is called
Kummer $(16)_{6}$-configuration.  On the other hand, 
${\rm Km}(J(C))$ is isomorphic to a complete intersection
of three quadrics in ${\bf P}^5$ on which 32 curves appear as 32 lines (Griffiths and Harris \cite{GH}).  

The authors would like to thank Professor I. Dolgachev for his useful comments.

\section{Preliminaries}
In this section, we recall some basic facts on the relation between
divisors and endomorphisms of a superspecial abelian surface.
Let $k$ be an algebraically closed field of characteristic $p > 0$
and $E$ be a supersingular elliptic curve defined over $k$.
We denote by $0$ the zero point of $E$.
We consider a superspecial abelian surface $A = E_{1} \times E_{2}$
with $E_{1} = E_{2} = E$. We set 
$X = E_{1}\times \{0\} +  \{0\}\times E_{2}$. 
$X$ is a principal polarization on $A$.
We set ${\cal O} = {\rm End}(E)$ and $B = {\rm End}(E)\otimes {\bf Q}$.
Then, $B$ is a quatenion division algebra over the rational number field ${\bf Q}$
with discriminant $p$, and ${\cal O}$ is a maximal order of $B$.
We denote by $\bar{a}$ the conjugate of $a \in B$.
For a divisor $L$, we have a homomorphism
$$
\begin{array}{cccc}
   \varphi_{L} :  &A  &\longrightarrow  &{\rm Pic}^{0}(A) \\
      & x & \mapsto & T_{x}^{*}L - L,
\end{array}
$$
where $T_{x}$ is the translation by $x \in A$.
We set 
$$
     H = \{
\left(
\begin{array}{cc}
\alpha  & \beta \\
\gamma & \delta
\end{array}
\right)
~\mid
~\alpha, \delta \in {\bf Z},~\gamma = \bar{\beta}
\}.
$$
Then we have the following theorem (cf. \cite{K1} and \cite{K2}).
\begin{th}\label{intersection}
The homomorphism
$$
\begin{array}{cccc}
j : & NS(A) & \longrightarrow  & H \\
   & L  &\mapsto  & \varphi_{X}^{-1}\circ\varphi_{L}
\end{array}
$$
is bijective.  By this correspondence, we have
$$
j(E_{1}\times \{0\}) = \left(
\begin{array}{cc}
0  & 0 \\
0 & 1
\end{array}
\right),~
j(\{0\}\times E_{2}) = \left(
\begin{array}{cc}
1  & 0 \\
0 &  0
\end{array}
\right).
$$
For $L_{1}, L_{2} \in NS(A)$ such that
$$
  j(L_{1}) =
\left(
\begin{array}{cc}
\alpha_{1}  & \beta_{1} \\
\gamma_{1} & \delta_{1}
\end{array}
\right),~
  j(L_{2}) =
\left(
\begin{array}{cc}
\alpha_{2}  & \beta_{2} \\
\gamma_{2} & \delta_{2}
\end{array}
\right),
$$
the intersection number  $(L_{1}, L_{2})$ is given by
$$
(L_{1}, L_{2}) =  \alpha_{1}\delta_{2} + \alpha_{2}\delta_{1}  - \gamma_{1}\beta_{2}
-\gamma_{2}\beta_{1}.
$$
In particular, for $L\in NS(A)$ such that 
$j(L) = \left(
\begin{array}{cc}
\alpha  & \beta \\
\gamma & \delta
\end{array}
\right)
$
we have
$$
\begin{array}{l}
L^{2} = 2\det \left(
\begin{array}{cc}
\alpha  & \beta \\
\gamma & \delta
\end{array}
\right)\\
(L, E_{1} ) = \alpha,~(L, E_{2}) = \delta .
\end{array}
$$
\end{th}
Let  $m : E\times E \rightarrow E$ be the addition of $E$, and we set
$$
\Delta = {\rm Ker}~ m.
$$
We have $\Delta = \{(P, - P)  ~\mid~P\in E\}$. For two endomorphisms
$a_{1}, a_{2} \in {\cal O} = {\rm End}(E)$, we set
$$
\Delta_{a_{1}, a_{2}} = (a_{1}\times a_{2})^{*}\Delta.
$$
Using this notation, we have $\Delta = \Delta_{1,1}$.
We have the following theorem (cf. \cite{K2})
\begin{th}\label{div}
$$
j(\Delta_{a_{1}, a_{2}}) =
\left(
\begin{array}{cc}
\bar{a}_{1}a_{1}  & \bar{a}_{1}a_{2}\\
\bar{a}_{2}a_{1} & \bar{a}_{2}a_{2}
\end{array}
\right)
$$
In particular, we have
$$
j(\Delta) =
\left(
\begin{array}{cc}
1  & 1\\
1 & 1
\end{array}
\right)
$$
\end{th}

\section{Supersingular elliptic curve in characteristic 2}
From here on, we consider algebraic varieties over an
algebraically closed field $k$ of characteristic 2, unless otherwise 
mentioned.  We summarize, in this section, facts 
on the supersingular elliptic curve in characteristic 2
which we will use later.  We also recall the facts on
the N\'eron-Severi group of supersingular K3 surface 
with Artin invariant 1
in charactersitic 2 (cf. Dolgachev and Kondo \cite{DK}) to explain
our aim.

We have, up to isomorphism, only one supersingular elliptic curve 
defined over $k$, which is given by the equation
$$
               y^{2} + y = x^{3}.
$$
We denote by $E$ a nonsingular complete model of the supersingular 
elliptic curve,
which is defined by
$$
  Y^{2}Z + YZ^{2} = X^{3},
$$
in the projective plane ${\bf P}^{2}$,
where $(X, Y, Z)$ is the homogeneous coordinate of ${\bf P}^{2}$.
In the affine model, let $(x_{1}, y_{1})$ and $(x_{2}, y_{2})$ be
two points on E. Then, the addition of E is given by
$$
\begin{array}{l}
x = x_{1} + x_{2} + ((y_{1} - y_{2})/(x_{1} - x_{2}))^{2},\\
y = ((y_{1} - y_{2})/(x_{1} - x_{2}))^{3}+
(x_{1}y_{1} + x_{2}y_{2})/(x_{1} -x_{2}) + 1.
\end{array}
$$
We denote by $[2]_{E}$ the multiplication by 2.  It is concretely
given by
$$
\begin{array}{l}
    x = x_{1}^{4},\\
    y = y_{1}^{4} + 1.
\end{array}
$$
We often use $2$ for $[2]_{E}$ for the sake of simplicity.
We denote by F (resp. V) the relative Frobenius morphism (resp.
the Vershiebung), which has the following relations:
$$
FV = 2,~V = \bar{F} = -F,~F^{2} = -2.
$$
We denote by $E({\bf F}_{q})$ the ${{\bf F}_{q}}$-rational points of $E$.
Then, we have
$$
\begin{array}{l}
E({\bf F}_{2}) = \{(0,1,0),(0,0,1), (0,1,1)\}, \\
E({\bf F}_{4}) = \{(0,1,0),(0,0,1), (0,1,1), (1, \omega,1),
 (1, \omega^{2},1), \\
\hspace*{2.0cm}(\omega, \omega,1), (\omega^{2}, \omega,1), 
(\omega, \omega^{2},1), (\omega^{2}, \omega^{2},1)\}.
\end{array}
$$
Here, $\omega$ is a primitive cube root of unity.
We set 
$$
P_{0} = (0,1,0), P_{1} = (0,0,1), P_{2} = (0,1,1).
$$
The point $P_{0}$ is the zero point of $E$, and the group of 
three torsion points of $E$ is given by
$$
 E({\bf F}_{4}).
$$

Next we recall some facts on the N\'eron-Severi group of supersingular K3 surfaces.
Let $S$ be a supersingular K3 surface with Artin invariant $\sigma$
in characteristic 2.  Note that the N\'eron-Severi group $NS{(S)}$ coincides 
with the Picard lattice for K3 surfaces.  
Together with the intersection form, $NS(S)$ is an even lattice of signature $(1,21)$ with the discriminant $-2^{2\sigma}$.  

In the following we assume that $S$ is a supersingular K3 surface with Artin invariant 1
in characteristic 2.   It is known that such K3 surface is unique up to isomorphisms
(cf. Rudakov and Shafarevich \cite{RS}, for instance).
To study the N\'eron-Severi group $NS(S)$, 
we use Conway's theorem on
the reflection group of an even unimodular lattice 
$L = {\rm II_{1,25}}$ of signature $(1,25)$.  We fix an orthogonal decomposition 
$L = \Lambda \perp U$,
where $\Lambda$ is the Leech lattice and $U$ is 
the hyperbolic plane.
Let $f, g$ be a standard basis of $U$,
i.e., $f^{2} = 0, g^{2} = 0$, and $\langle f, g\rangle = 1$,
and we denote each vector $x \in L$ by $(\lambda, m, n)$,  where
$x = \lambda + mf + ng$ with $\lambda \in \Lambda$. 
The vector $\rho = (0, 0, 1)$ is
called the Weyl vector. We have $\rho^{2} = 0$.
For $r \in L$ with $r^{2} = -2$, $r$ is called a Leech root
if $\langle \rho, r \rangle = 1$. We denote by $\Delta (L)$
the set of Leech roots of $L$. 
We denote by $W(L)$ the group generated by the reflections 
in the orthogonal group
${\rm O}(L)$. We set
$$
P(L) = \{x \in L \otimes {\bf R} ~\mid~ x^{2}>0\}.
$$
Then, $P(L)$ has two connected components, and we denote by
$P(L)^{+}$ one of them.  Then, $W(L)$ acts on $P(L)^{+}$. We set 
$$
C = \{x \in P(L)^{+}~\mid~ \langle x, r\rangle > 0
~\mbox{for any}~r \in  \Delta (L)\}.
$$
Leech roots are important because of the following theorem.
\begin{th}[Conway \cite{Co}]
The closure $\bar{C}$ of $C$ in $P(L)^+$ is a fundamental domain of $W(L)$.
\end{th}
There exists a primitive embedding of the N\'eron-Severi lattice $NS{(S)}$ into $L$
such that $P^+(L)$ contains the positive cone $P^+(S)$ of $S$.  Under the embedding, 
the restriction $C \cap P^+(S)$ is a finite polyhedron whose faces consist of 42 hyperplanes
defined by 42 $(-2)$-vectors in $NS(S)$ and 168 hyperplanes defined by 168 $(-4)$-vectors in $NS(S)$ (see \cite{B}).  
Each of these 42 $(-2)$-vectors can be represented by a smooth rational curve on $S$, 
and these 42 curves are divided into two families ${\cal A}$, ${\cal B}$ 
as mentioned in the introduction.  The group of symmetries of the 42 $(-2)$-vectors
is isomorphic to ${\rm PGL}(3, {\bf F}_{4}) \cdot {\bf Z}/2{\bf Z}$ where 
${\rm PGL}(3, {\bf F}_{4})$
preserves each family and ${\bf Z}/2{\bf Z}$ switches $\cal A$ and $\cal B$.
The group  ${\rm PGL}(3, {\bf F}_{4}) \cdot {\bf Z}/2{\bf Z}$ can be represented by 
automorphisms of $S$.  The projection of the Weyl vector $\rho$ into 
$NS(S)\otimes {\bf Q}$ sits in $NS(S)$
and is an ample class of degree 14.  Obviously, with respect to this ample class, 
42 curves are lines.  For more details, we refer the reader to \cite{DK}.

\section{Construction of elliptic curves on an abelian surface}
We use the notation in Section 3.  Let $E$ be the supersingular elliptic curve.  
Then $E$ has the following automorphism $\sigma$ of order 3:
$$  \sigma : x\mapsto \omega x,~y \mapsto y.$$
We consider a superspecial abelian surface
$$A = E \times E.$$
Then, the automorphism $\sigma \times \sigma^{2}$ acts on $A$, and the minimal resolution ${\rm GKm}(A)$ of the quotient surface $A/\langle \sigma \times \sigma^{2}\rangle$ is a K3 surface.
In this section, we construct 24 elliptic curves on $A$ which are invariant under the automorphism $\sigma \times \sigma^{2}$.
In the next section, we will see that these elliptic curves give 24 nonsingular rational curves on ${\rm GKm}(A)$ which correspond to 
a part of 42 smooth rational curves which correspond with Leech roots.

The translation by a 3-torsion point $P_{1}$ of $E$ is given by
$$\tau: ~x \mapsto y/x^{2},~y \mapsto y/x^{3}.$$
It acts on the function field $k(E) = k(x, y)$ of $E$ and the invariant field $k(E)^{\langle \tau \rangle}$ is given by
$$k(E)^{\langle \tau \rangle} = k((x^{3} + 1)/x^{2},y + 1 + (1/x^{3}) + \omega).$$
We set 
$$\begin{array}{l}w =(x^{3} + 1)/x^{2}, ~ z = y + 1 + (1/x^{3}) + \omega.\end{array}$$
Then, we know that the quotient curve $E/\langle \tau \rangle$ is defined by
$$z^{2} + z = w^{3},$$
which coincides with the original elliptic curve $E$.  We denote by $\pi$ the projection given by
$$
(x, y) \mapsto (w, z).
$$
We also consider the following automorphism of E:
$$
    \theta : x\mapsto x + 1,~ y \mapsto y + x + \omega.
$$
Then, we have $\theta^{2} = -{\rm  id}_{E}$, and by direct calculation we have
$$
\begin{array}{l}
         \bar{\sigma} = \sigma^{2}\\
         F\circ \sigma = \sigma^{2}\circ F \\
         \pi = \theta\circ ({\rm id}_{E} -F) = 2 \sigma + 1\\
         \bar{\pi} = -\pi\\
         F\circ \pi = -\pi \circ F \\
         F = \sigma\circ \theta  - \theta\circ \sigma  = 
\theta\circ \sigma^{2} - \sigma^{2}\circ \theta\\
        {\rm id}_{E} = \theta\circ \sigma - \sigma^{2}\circ \theta.
\end{array}
$$

We consider eight homomorphisms from $E$ to $A$ which are given as follows:
$$
\begin{array}{l}(1) \quad (0, {\rm id}_{E}), \quad (2) \quad (F, \sigma), \quad    (3) \quad (V, \sigma^{2}), \quad (4) \quad (V, \pi), \\(1') \quad ({\rm id}_{E}, 0),\quad (2')\quad (\sigma^{2}, F), \quad (3') \quad (\sigma, V), \quad (4') \quad (\pi, F).\end{array}
$$
All these homomorphisms are injective and the images of them are elliptic curves which are invariant under the action of $\sigma \times \sigma^{2}$. 
 We denote the elliptic curves respectively by
$$
\begin{array}{l}(1)~ E_{0}, \quad (2)~ F_{0}, \quad (3)~ V_{0}, \quad (4)~ \pi_{0}, \\  (1')~ E_{0}^{\prime},\quad (2')~ F_{0}^{\prime}, \quad (3')~ V_{0}^{\prime}, \quad  (4')~ \pi_{0}^{\prime}.
\end{array}
$$
Using the notation in Section 2, we have
$$
\begin{array}{l}
(1)~ E_{0} = \{0\}\times E_{2}, \quad (2)~ F_{0}= \Delta_{\sigma^{2}, V}, 
\quad (3)~ V_{0} =\Delta_{\sigma, F}, \quad (4)~ \pi_{0}= \Delta_{-\pi, F}, \\  
(1')~ E_{0}^{\prime}= E_{1}\times \{0\},
\quad (2')~ F_{0}^{\prime}= \Delta_{V, \sigma}, 
\quad (3')~ V_{0}^{\prime}=\Delta_{F, \sigma^{2}}, 
\quad  (4')~ \pi_{0}^{\prime}=\Delta_{V, -\pi}.
\end{array}
$$
and
$$
\begin{array}{ll}
(1)~ j(E_{0}) = 
\left(
\begin{array}{cc}
1  & 0 \\
0 & 0
\end{array}
\right),&
(2)~ j(F_{0})= \left(
\begin{array}{cc}
1  & \sigma V \\
F\sigma^{2} & 2
\end{array}
\right),\\
 (3)~j(V_{0}) =\left(
\begin{array}{cc}
1  & \sigma^{2}F \\
V \sigma & 2
\end{array}
\right), &
(4)~ j(\pi_{0})= \left(
\begin{array}{cc}
3  & \pi F \\
F\pi & 2
\end{array}
\right), \\  
(1')~ j(E_{0}^{\prime})=\left(
\begin{array}{cc}
0  & 0 \\
0 & 1
\end{array}
\right), &
(2')~ j(F_{0}^{\prime})= \left(
\begin{array}{cc}
2  & F\sigma\\
\sigma^{2}V & 1
\end{array}
\right),\\
(3')~ j(V_{0}^{\prime})=\left(
\begin{array}{cc}
2  & V\sigma^{2} \\
\sigma F & 1
\end{array}
\right),& 
(4')~ j(\pi_{0}^{\prime})=\left(
\begin{array}{cc}
2  & V\pi \\
\pi V & 3
\end{array}
\right).
\end{array}
$$
Therefore, using the intersection theory in Section 2, we have intersection numbers
of these elliptic curves as follows:

\begin{center}
\begin{tabular}{| l |c|c|c|c|c|c|c|c|}
\hline
     & $E_{0}$ & $F_{0}$ & $V_{0}$ & $\pi_{0}$ & $E_{0}^{\prime}$ & $F_{0}^{\prime}$ & $V_{0}^{\prime}$ & $\pi_{0}^{\prime}$\\
\hline
$E_{0}$ & 0 & 2 & 2& 2& 1& 1 & 1 & 3\\
\hline
$F_{0}$ & 2 & 0 & 2 & 2 & 1 & 3 & 1 & 2\\
\hline
$V_{0}$ & 2 & 2& 0 & 2& 1 & 1 & 3 & 2\\
\hline
$\pi_{0}$ & 2 & 2& 2 & 0 & 3 & 1 & 1 & 1\\
\hline
$E_{0}^{\prime}$ &  1 & 1& 1& 3& 0 & 2 & 2 & 2\\
\hline
$F_{0}^{\prime}$ & 1 & 3 & 1& 1& 2 & 0 & 2 & 1\\
\hline
$V_{0}^{\prime}$ & 1 & 1 & 3& 1& 2 & 2 & 0 & 1\\
\hline
$\pi_{0}^{\prime}$ & 3 & 2 & 2& 1& 2 & 1 & 1 & 0\\
\hline
\end{tabular}
\end{center}

We say that the elliptic curves $E_{0}, F_{0}, V_{0}, \pi_{0}$ (resp. $E_{0}^{\prime}, F_{0}^{\prime}, V_{0}^{\prime}, \pi_{0}^{\prime}$) 
are the elliptic curves in the first group (resp. the elliptic curves in the second group).  
Note that the elliptic curves in the first group intersect the elliptic curves in the second group with multiplicity 1
at each intersection point and that the elliptic curves in the first group (resp. the elliptic curves in the second group) 
intersect each other with multiplicity 2.  On $A$, $\sigma \times \sigma^{2}$ has 9 fixed points which are given 
by $P_{i}\times P_{j}$ $(i, j = 0, 1, 2)$.  We set $P_{ij} = P_{i}\times P_{j}$ ($i, j = 0, 1, 2$). 
These points coincide with the ${\bf F}_{2}$-rational points $A({\bf F}_{2})$ of $A$.  
We denote by $T_{ij}$ the translation by $P_{ij}$, and we set
$$
\begin{array}{l}(1) ~E_{i} =T_{i0}E_{0}, \quad (2) ~F_{i} =T_{i0}F_{0}, \quad(3) ~V_{i} =T_{i0}V_{0}, \quad (4)~\pi_{i}=T_{0i}\pi_{0},\\  (1')~ E_{j}^{\prime}=T_{0j}E_{0}^{\prime},\quad (2')~ F_{j}^{\prime}=T_{0j}F_{0}^{\prime},\quad (3')~ V_{j}^{\prime}=T_{0j}V_{0}^{\prime},\quad(4')~ \pi_{j}^{\prime}=T_{j0}\pi_{0}^{\prime}, \end{array}
$$
with $i= 1, 2$ and $j = 1, 2$.   All these curves are invariant under the action of $\sigma \times \sigma^{2}$.  Thus, we have in total  
24 elliptic curves which are all invariant under the action of $\sigma \times \sigma^{2}$. 
Among these 24 curves, there are 8 elliptic curves which pass through the point $P_{ij}$, 
and each of these elliptic curves passes through 3 points among points in $A({\bf F}_{2}) = \{ P_{ij}~\mid~i, j = 0, 1, 2\}$. 
We give here the list of elliptic curves which pass through the point 
$P_{ij} =P_{i}\times P_{j}$ $(1\leq i, j \leq 3)$.
\begin{center}
\begin{tabular}{| l |c|c|c|}
\hline
     & $P_{0}$ & $P_{1}$ & $P_{2}$ \\
\hline
$P_{0}$ & 
$
\begin{array}{c}
E_{0}, F_{0}, V_{0}, \pi_{0}\\
        E_{0}', F_{0}', V_{0}', \pi_{0}'   
\end{array}   
$       
 &
$
 \begin{array}{c}
            E_{1}, F_{1}, V_{1}, \pi_{0}\\
        E_{0}', F_{2}', V_{1}', \pi_{1}'     
\end{array}
$   
 &
$
\begin{array}{c}
E_{2}, F_{2}, V_{2}, \pi_{0}\\
        E_{0}', F_{1}', V_{1}', \pi_{2}'   
\end{array}   
$ \\
\hline
$P_{1}$ & 
$
\begin{array}{c}
E_{0}, F_{2}, V_{1}, \pi_{1}\\
        E_{1}', F_{1}', V_{1}', \pi_{0}'   
\end{array}   
$       
 &
$
 \begin{array}{c}
E_{1}, F_{0}, V_{2}, \pi_{1}\\
        E_{1}', F_{0}', V_{2}', \pi_{1}'   
\end{array}   
$
 &
$
\begin{array}{c}
E_{2}, F_{1}, V_{0}, \pi_{1}\\
        E_{1}', F_{2}', V_{0}', \pi_{2}'   
\end{array}
$   
\\

\hline
$P_{2}$ & 
$
\begin{array}{c}
E_{0}, F_{1}, V_{2}, \pi_{2}\\
        E_{2}', F_{2}', V_{2}', \pi_{0}'   
\end{array} 
$         
 & 
$
\begin{array}{c}
E_{1}, F_{2}, V_{0}, \pi_{2}\\
        E_{2}', F_{1}', V_{0}', \pi_{1}'   
\end{array} 
$  
 &
$
\begin{array}{c}
E_{2}, F_{0}, V_{1}, \pi_{2}\\
        E_{2}', F_{0}', V_{1}', \pi_{2}'   
\end{array}   
$\\
\hline
\end{tabular}
\end{center}
For example, we see that the elliptic curve $F_{1}$ passes through 3-torsion points
$P_{1}\times P_{0}$, $P_{2}\times P_{1}$ and $P_{0}\times P_{2}$ from this table.

Now we say that elliptic curves with $i = 1, 2$ (resp. $j = 1, 2$) are also elliptic curves in the first group (resp. in the second group).  
Since the intersection numbers don't change by any translation,  
we know that the elliptic curves in the first group intersect the elliptic curves in the second group with multiplicity 1 
at each intersection point and that the elliptic curves in the first group (resp. the elliptic curves in the second group) 
intersect each other with multiplicity 2 at each intersection point. 
Besides three points in $A({\bf F}_{2})$, each of these 24 elliptic curves contains 6 more points of order 3 of $A$, 
which are contained in $A({\bf F}_{4})$. The group $\langle \sigma \times \sigma^{2}\rangle$ of order 3 acts on 
these six points and each set of the six points has two orbits. 
The group $\langle \sigma \times \sigma^{2}\rangle$ also acts on the intersection points of each two elliptic curves 
which we constructed. For example, $F_{0}$ intersects $F_{0}^{\prime}$ at the points $P_{00}, P_{11}, P_{22}$. 
and we have $(F_{0}\cdot F_{0}^{\prime}) = 3$ as in the table above. 
Therefore, we have 
$$
 (F_{2}\cdot F_{0}^{\prime}) = (F_{0}\cdot F_{0}^{\prime})= 3.
$$
By a direct calculation, we see
$$
F_{2}\cap F_{0}^{\prime} = \{(\omega^{2},\omega^{2})\times (1, \omega),(1,\omega^{2})\times (\omega^{2}, \omega), (\omega,\omega^{2})\times (\omega, \omega)\}.
$$
The group $\langle \sigma \times \sigma^{2}\rangle$ of order 3 acts on these three points transitively.  
The elliptic curve $\pi_{0}$ intersects $E_{0}^{\prime}$ at three points $P_{00}, P_{10}, P_{20}$. 
These are fixed points of the group $\langle \sigma \times \sigma^{2}\rangle$. 
The elliptic curve $\pi_{0}$ intersects  $E_{1}^{\prime}$ (resp. $E_{2}^{\prime}$) at three 3-points $(1,\omega)\times P_{1}$, 
$(\omega^{2}, \omega)\times P_{1}$,  $(\omega, \omega)\times P_{1}$ 
(resp. three 3-torsion points $(1,\omega^{2})\times P_{2}$, $(\omega^{2}, \omega^{2})\times P_{2}$, $(\omega, \omega^{2})\times P_{2})$. 
The group $\langle \sigma \times \sigma^{2}\rangle$ of order 3 acts on these three points transitively.

We give the list of elliptic curves that shows which 3-torsion points 
they pass through. Snice we already gave the list of elliptic curves which pass
through $P_{i}\times P_{j}$ ($1\leq i, j \leq 3$) defined over ${\bf F}_{2}$,
we give the list for the points in $A({\bf F}_{4})\setminus A({\bf F}_{2})$.

{\scriptsize
\begin{center}
\begin{tabular}{| l |c|c|c|c|c|c|c|c|c|}
\hline
     & $(1, \omega)$ & $(1, \omega^2)$ & $(\omega, \omega)$ & $(\omega, \omega^2)$ & $(\omega^2, \omega)$ & $(\omega^2, \omega^2)$ & $P_{0}$ & $P_{1}$ & $P_{2}$\\
\hline
$(1, \omega)$ &$V_{1},V_{2}'$ &$F_{1},F_{1}'$ & $V_{2},V_{0}'$& $F_{0},F_{2}'$& 
$V_{0},V_{1}'$& $F_{2},F_{0}'$ & $E_{0},\pi_{2}'$ & $E_{1},\pi_{0}'$&$E_{2},\pi_{1}'$\\
\hline
$(1, \omega^2)$ &$F_{2},F_{2}'$ &$V_{2},V_{1}'$ & $F_{0},F_{1}'$& $V_{1},V_{0}'$& 
$F_{1},F_{0}'$& $V_{0},V_{2}'$ & $E_{0},\pi_{1}'$ & $E_{1},\pi_{2}'$&$E_{2},\pi_{0}'$\\
\hline
$(\omega, \omega)$ &$V_{2},V_{0}'$ &$F_{0},F_{2}'$ & $V_{0},V_{1}'$& $F_{2},F_{0}'$& 
$V_{1},V_{2}'$& $F_{1},F_{1}'$ & $E_{0},\pi_{2}'$ & $E_{1},\pi_{0}'$&$E_{2},\pi_{1}'$\\
\hline
$(\omega, \omega^2)$&$F_{0},F_{1}'$ &$V_{1},V_{0}'$ & $F_{1},F_{0}'$& $V_{0},V_{2}'$& 
$F_{2},F_{2}'$& $V_{2},V_{1}'$ & $E_{0},\pi_{1}'$ & $E_{1},\pi_{2}'$&$E_{2},\pi_{0}'$\\
\hline
$(\omega^2, \omega)$ &$V_{0},V_{1}'$ &$F_{2},F_{0}'$ & $V_{1},V_{2}'$& $F_{1},F_{1}'$& 
$V_{2},V_{0}'$& $F_{0},F_{2}'$ & $E_{0},\pi_{2}'$ & $E_{1},\pi_{0}'$&$E_{2},\pi_{1}'$\\
\hline
$(\omega^2, \omega^2)$ &$F_{1},F_{0}'$ &$V_{0},V_{2}'$ & $F_{2},F_{2}'$& $V_{2},V_{1}'$& 
$F_{0},F_{1}'$& $V_{1},V_{0}'$ & $E_{0},\pi_{1}'$ & $E_{1},\pi_{2}'$&$E_{2},\pi_{0}'$\\
\hline
$P_{0}$ &$\pi_{2}, E_{0}'$ &$\pi_{1}, E_{0}'$& $\pi_{2}, E_{0}'$& $\pi_{1}, E_{0}'$& 
$\pi_{2}, E_{0}'$& $\pi_{1}, E_{0}'$ &  & &\\
\hline
$P_{1}$ &$\pi_{0}, E_{1}'$ &$\pi_{2}, E_{1}'$& $\pi_{0}, E_{1}'$& $\pi_{2}, E_{1}'$& 
$\pi_{0}, E_{1}'$& $\pi_{2}, E_{0}'$ &  & &\\
\hline
$P_{2}$ &$\pi_{1}, E_{2}'$ &$\pi_{0}, E_{2}'$& $\pi_{1}, E_{2}'$& $\pi_{0}, E_{2}'$& 
$\pi_{1}, E_{2}'$& $\pi_{0}, E_{2}'$ &  & &\\
\hline
\end{tabular}
\end{center}
}
For example, at the 3-torsion point $(1, \omega^{2})\times (1, \omega)$
we see from this table that the elliptic curve $F_{1}$ intersects $F_{1}'$ .
We also see that for each point of $A({\bf F}_{4})\setminus A({\bf F}_{2})$,
there exist just two elliptic curves passing through the point 
among these 24 elliptic curves.

\section{A generalized Kummer surface}
We use the notation in Sections 3 and 4.  We consider the superspecial abelian surface $A = E\times E$  and 
an automorphism $\sigma \times \sigma^{2}$ of $A$. The automorphism $\sigma \times \sigma^{2}$ has nine fixed points, 
which are given by $A({\bf F}_{2})$. Let $A/\langle \sigma \times \sigma^{2}\rangle$ be the quotient surface of $A$ 
by the group of order 3 which is generated by $\sigma \times \sigma^{2}$, and let ${\rm GKm}(A)$ 
be the minimal resolution of $A/\langle \sigma \times \sigma^{2}\rangle$.  
%Then, ${\rm GKm}(A)$ is a supersingular K3 surface with Artin invariant 1 (see Schr\"oer \cite{Sch}, 
%Dolgachev and Kondo \cite{DK} and Proposition 5.1 in this paper).  
Since $A/\langle \sigma \times \sigma^{2}\rangle$
has 9 rational double points of type $A_{2}$, on ${\rm GKm}(A)$ there exist 18 exceptional curves, 
which are nonsingular rational curves.  We take the images of 24 elliptic curves which we constructed in the previous section.  
They are also nonsingular rational curves. Therefore, we have in total 42 nonsingular rational curves on ${\rm GKm}(A)$.  
To examine these 42 rational curves, first we blow up at nine points $A({\bf F}_{2})$ of A.  
We denote by $A^{\prime}$ the surface obtained from $A$.  Then, the automorphism $\sigma \times \sigma^{2}$ 
can be lifted to an automorphism of $A^{\prime}$.  It has two fixed points on each exceptional curves.  
We blow up again these 18 fixed points on $A^{\prime}$. We denote by $A^{\prime \prime}$ the surface obtained from $A^{\prime}$.  Then, the automorphism can be again lifted to an automorphism of $A^{\prime \prime}$.  
We denote by $\eta$ the automorphism on $A^{\prime \prime}$ which we got in this way.  
Then $\eta$ has 18 fixed exceptional curves which we got in the second step.  
We set ${\tilde A} = A^{\prime \prime}/\langle \eta \rangle$.  
Then, ${\tilde A}$ is nonsingular and contains 9 exceptional curves of the first kind 
which derive from the exceptional curves on $A^{\prime}$. 
We blow down the 9 exceptional curves. Then, we have the generalized Kummer surface ${\rm GKm}(A)$ 
which we already got above.  We summarize the above procedure as in the following diagram.
$$
\begin{array}{ccccc}
\  A^{\prime\prime}& \stackrel{f}{\longrightarrow} &  A^{\prime}& \stackrel{f'}{\longrightarrow} & A\\   \quad  \downarrow{g}    &                 &             &     & \downarrow \\  {\tilde A}  & \stackrel{h}{\longrightarrow} & {\rm GKm}(A) & \longrightarrow &              A/\langle \sigma \times \sigma^{2}\rangle.   
\end{array}
$$
On $A^{\prime\prime}$, we take proper transforms of 24 elliptic curves which we constructed in Section 3. 
They are divided into two groups, the curves in the first group and the curves in the second group, 
as in the previous section.  We consider the images of these curves by the morphism $h\circ g$.  
We also take the images of 18 exceptional curves by $h\circ g$ which are fixed by the automorphism $\eta$.  
All curves of the images are nonsingular and rational.  
The images of exceptional curves are divided into two groups, that is, 
the curves which do not intersect the curves in the first group and the curves which do not intersect the curves in 
the second group.
Therefore, we add the former images of exceptional curves (resp.  the latter images of exceptional curves) into the the first group 
(resp. into the second group).  Then, we have the following theorem.

\begin{th}
The $42$ curves which we constructed above on ${\rm GKm}(A)$ make two families $\cal{A}$, $\cal{B}$,each consisting of $21$ disjoint smooth rational curves. Each memberin one family meets exactly five members in another family.
\end{th}
\proof{
By the blowings-up $f'$, the proper transforms of the elliptic curves in the first group (resp. in the second group) 
passing through the point $P_{ij}$ become 4 elliptic curves which intersect at the same point 
(resp. at the same point which is different from the point for the first group).  
The blowings-up $f$ are the ones whose centers are the 18 intersection points of the elliptic curves.  
On $A^{\prime \prime}$, elliptic curves in the first group (resp. elliptic curves in the second group) 
are disjoint from each other.  We consider the elliptic curves and 27 exceptional curves 
of the morphism $f'\circ f$. The images by the projection $g$ of 9 exceptional curves on $A^{\prime \prime}$ 
which do not intersect the elliptic curves are exceptional curves of the first kind.  
By the blowings-down $h$ those exceptional curves collapse.  On $A^{\prime \prime}$, each elliptic curve in the first group 
meets two elliptic curves in the second group respectively on three points which derive from 3-torsion points on $A$ 
which are not defined over ${\bf F}_{2}$.  Each of these 3 points comes to a point on ${\tilde A}$. 
The situation is same for the elliptic curves in the second group. Therefore, on ${\tilde A}$ the image by $g$ of 
each elliptic curve in the first group (resp. in the second group) intersects two images of elliptic curves in the second group 
(resp. in the first group) and intersects three images of exceptional curves in the second group (resp. in the first group).  
Therefore, in total it intersects 5 curves in the second group (resp. in the first group) transversely.  
By the blowings-down $h$, the images of the exceptional curves in the first group (resp. in the second group) still intersect 
5 curves in the second group (resp. in the first group) transversely.}

The following corollary is well-known (cf. Shimada \cite{S}).
\begin{cor} 
{\it There exists a normal $K3$ surface with $21$ rational double points of type $A_{1}$
in characteristic two. The number $21$ is the maximal one of isolated singular points on $K3$ surface.}
\end{cor}
\proof{
We can collapse 21 rational curves in the first group.  The singular points are clearly rational double points of type $A_{1}$.
By Hodge index theorem, we know that 21 is the maximal number of isolated singular points of K3 surface.
}
Next we summarize the intersection numbers of 42 curves which we have constructed.  
We consider the elliptic curves on $A$ which we constructed in Section 2,  and we use the same notation 
as the original elliptic curves for the rational curves on ${\rm GKm}(A)$ which derive from the elliptic curves. 
Then, on ${\rm GKm}(A)$, the intersection numbers for the rational curves are given by
$$
\begin{array}{l}
F_{i}\cdot F_{j}^{\prime} = 1 -\delta_{ij},\\
V_{i}\cdot V_{j}^{\prime} = 1 - \delta_{ij}, \\
\pi_{i}\cdot E_{j}^{\prime} = 1 - \delta_{ij},\\
E_{i} \cdot \pi_{j}^{\prime} = 1 - \delta_{ij}
\end{array}
$$
($i, j = 0, 1, 2$) and for others from the elliptic curves, the intersection numbers are 0.  
Here, $\delta_{ij}$ is Kronecker's delta.  We denote by $\ell_{qr}$ (resp. $\ell_{qr}^{\prime}$) ($q, r = 0, 1, 2$) 
the images of exceptional curves in the first group (resp. in the second group) 
which derive from the blowings-ups over the point $P_{qr}$.  Then, $E_{i},F_{i}, V_{i}, \pi_{i}$ ($i = 0, 1, 2$) 
(resp. $E_{j}^{\prime},F_{j}^{\prime}, V_{j}^{\prime}, \pi_{j}^{\prime}$ ($j = 0, 1, 2$)) and $\ell_{qr}$ ($q,r = 0, 1, 2$) 
(resp. $\ell_{qr}^{\prime}$ ($q, r = 0, 1, 2$)) are the disjoint rational curves in the first group (resp. in the second group). 
The intersection numbers of $\ell_{qr}, \ell_{qr}^{\prime}$ are given as follows:
$$
\begin{array}{l}\ell_{qr}\cdot \ell_{st}^{\prime} = \delta_{qs}\delta_{rt},\\
\ell_{qr}\cdot E_{i}^{\prime} = \delta_{ri},~\ell_{qr}\cdot \pi_{i}^{\prime} = \delta_{qi},~\ell_{qr}^{\prime}\cdot E_{i} = \delta_{qi},~\ell_{qr}^{\prime}\cdot \pi_{i} = \delta_{ri},\\\ell_{qr}\cdot F_{0}^{\prime} =\delta_{qr},~\ell_{qr}^{\prime}\cdot F_{0} =\delta_{qr},\\
\ell_{qr}\cdot F_{1}^{\prime} =1~\mbox{for}~(q, r) = (0,1),(1,2),(2,0),~ \ell_{qr}\cdot F_{1}^{\prime} = 0 ~\mbox{for other}~q,r,\\
\ell_{qr}^{\prime}\cdot F_{1}= 1 ~\mbox{for}~(q, r) = (0,2),(1,0),(2,1),~\ell_{qr}^{\prime}\cdot F_{1}= 0 ~\mbox{for other}~q,r,\\
\ell_{qr}\cdot F_{2}^{\prime}= 1 ~\mbox{for}~(q, r) = (0,2),(1,0),(2,1),~\ell_{qr}^\cdot F_{2}^{\prime} = 0 ~\mbox{for other}~q,r,\\
\ell_{qr}^{\prime}\cdot F_{2}= 1 ~\mbox{for}~(q, r) = (0,1),(1,2),(2,0),~\ell_{qr}^{\prime}\cdot F_{2} = 0 ~\mbox{for other}~q,r,\\
\ell_{qr}\cdot V_{0}^{\prime} =1~\mbox{for}~(q, r) = (0,0),(1,2),(2,1),~ \ell_{qr}\cdot V_{0}^{\prime} = 0 ~\mbox{for other}~q,r,\\\ell_{qr}^{\prime}\cdot V_{0}= 1 ~\mbox{for}~(q, r) = (0,0),(1,2),(2,1),~\ell_{qr}^{\prime}\cdot V_{0}= 0 ~\mbox{for other}~q,r,\\\ell_{qr}\cdot V_{1}^{\prime}= 1 ~\mbox{for}~(q, r) = (0,1),(1,0),(2,2),~\ell_{qr}\cdot V_{1}^{\prime} = 0 ~\mbox{for other}~q,r,\\\ell_{qr}^{\prime}\cdot V_{1}= 1 ~\mbox{for}~(q, r) = (0,1),(1,0),(2,2),~\ell_{qr}^{\prime}\cdot V_{1} = 0 ~\mbox{for other}~q,r,\\\ell_{qr}\cdot V_{2}^{\prime} =\delta_{(q+r),2},~\ell_{qr}^{\prime}\cdot V_{2} =\delta_{(q+r),2}.\\\end{array}
$$
Since the self-intersection number of the divisor
$$F_{0} + F_{1}^{\prime}+ F_{2} + F_{0}^{\prime}+ F_{1}+ F_{2}^{\prime}$$
is zero, the linear system associated with this divisor gives a structure of elliptic surface on 
${\rm GKm}(A)$.  Besides this divisor, this elliptic surface has 3 singular fibers:
$$\begin{array}{l}V_{0} + V_{1}^{\prime}+ V_{2} + V_{0}^{\prime}+ V_{1}+ V_{2}^{\prime},\\\pi_{0} + E_{1}^{\prime}+ \pi_{2} + E_{0}^{\prime}+ \pi_{1}+ E_{2}^{\prime},\\E_{0} + \pi_{1}^{\prime}+ E_{2} + \pi_{0}^{\prime}+ E_{1}+ \pi_{2}^{\prime}  .\end{array}$$
All these singular fibers are of type ${\rm I}_{6}$ (cf. Kodaira\cite{Ko}).  Since ${\rm GKm}(A)$ is a K3 surface, 
the second Chern number $c_{2}({\rm GKm}(A))$ is equal to 24. 
Therefore, considering the discriminant of this elliptic surface, we conclude that we have no more singular fibers.  
The 18 rational curves $\ell_{qr}$ and $\ell_{qr}^{\prime}$ ($q, r = 0, 1, 2$) give sections of this elliptic surface.  
Two sections $\ell_{qr}$ and $\ell_{qr}^{\prime}$ intersect each other with intersection number 1 
on the only one fiber which is a supersingular elliptic curve.  These 9 points are all of intersection points 
for these 18 sections.  These 18 sections generate the Modell-Weil group which is isomorphic to ${\bf Z}/6{\bf Z} \times {\bf Z}/3{\bf Z}$.
Since the elliptic fibration has four singular fibers of type ${\rm I}_6$ and a section, ${\rm GKm}(A)$ has the Picard number 22, and hence
is supersingular.  By the Tate-Shioda formula, the discriminant of the Picard lattice is $-6^4/18^2 = -2^2$, and hence ${\rm GKm}(A)$ has
Artin invariant 1.  Thus we have proved the following theorem (see also \cite{Sch}).

\begin{th}
The generalized Kummer surface ${\rm GKm}(A)$ is a supersingular $K3$ surface with Artin invariant $1$.
\end{th}

\noindent
By the uniqueness of the supersingular K3 surface with Artin invariant 1, ${\rm GKm}(A)$ is isomorphic to $S$ mentioned in the
Introduction.

Finally we shall give a few remarks.
Consider the pencil of cubics
$$
   X^{3} + Y^{3} + Z^{3} + \mu XYZ = 0.
$$
After blowing up the base points, we have a rational elliptic surface $V$ with four singular fibers of type ${\rm I}_3$.
In \cite{DK}, \S 4.3, Dolgachev and the second author showed that by taking the inseparable double cover of the base of the elliptic fibration $V$, 
one can obtain the supersingular K3 surface $S$ with the Artin invariant 1 which has a structure of an elliptic fibration with four singular fibers
of type ${\rm I}_6$ and 18 sections.  Also Ito {\rm \cite{It}} gave an equation 
of an elliptic K3 surface with four singular fibers of type ${\rm I}_6$  and 18 sections which is obtained by the inseparable base change
of a rational elliptic surface with four singular fibers of type ${\rm I}_3$.

On the other hand,
as mentioned in Introduction, $S$ is isomorphic to the surface
in ${\bf P}^2\times {\bf P}^2$ defined by the equations of bidegree $(1,2)$ and $(2,1)$.
Consider the involution of ${\bf P}^2\times {\bf P}^2$ defined by
$$((x_0,x_1,x_2), (y_0,y_1,y_2)) \to ((y_0, y_1, y_2), (x_0,x_1,x_2)),$$
which induces an involution $\iota$ of $S$.  Note that $\iota$ switches the set $\cal{A}$
and $\cal{B}$, and is a generator of ${\bf Z}/2{\bf Z}$ in the group
${\rm PGL}(4, {\bf F}_{2^{2}}) \cdot {\bf Z}/2{\bf Z}$ (see Section 2).
The set of fixed points of $\iota$ is given by
$$E : x_0^3 +x_1^3+x_2^3 = 0$$
which is a supersingular elliptic curve in characteristic $2$.  
The elliptic fibration given by the linear system $|E|$ has four singular fibers of type ${\rm I}_6$
and $18$ sections.  
We remark that 24 components and 18 sections correspond to 
the images of 24 elliptic curves and 18 exceptional curves constructed in 
Sections 4 and 5.
The dual graph of 18 sections is isomorphic to $A_2^{\oplus 9}$.  Each pairs of sections forming
$A_2$ meet at one of 9 inflection points on $E$.
The quotient $S/\langle \iota \rangle$ is a rational elliptic surface with
four singular fibers of type ${\rm I}_3$ and $9$ sections.  
The double cover $S \to S/\langle \iota \rangle$ is separable contrary to the above inseparable construction of $S$ from $V$
due to Dolgachev and  Kondo \cite{DK}.

\vspace{0.5cm}
\noindent
T.\ Katsura: Faculty of Science and Engineering, Hosei University,
Koganei-shi, Tokyo 184-8584, Japan

\noindent
E-mail address: toshiyuki.katsura.tk@hosei.ac.jp

\vspace{0.3cm}
\noindent
S.\ Kond\=o: Graduate School of Mathematics, Nagoya University, Nagoya 464-8602, Japan

\noindent
E-mail address: kondo@math.nagoya-u.ac.jp

\end{document}